\documentclass{amsart}
\usepackage{amssymb}
\usepackage{amsfonts}
\usepackage{amsbsy}

\usepackage{latexsym}
\usepackage[mathscr]{eucal}
\usepackage{verbatim}

\theoremstyle{plain}
\newtheorem{theorem}{Theorem}[section]
\newtheorem{corollary}[theorem]{Corollary}
\newtheorem{lemma}[theorem]{Lemma}
\newtheorem{proposition}[theorem]{Proposition}
\newtheorem{remark}[theorem]{Remark}
\newtheorem{definition}[theorem]{Definition}

\newcommand{\be}{\begin{equation}\label}
\newcommand{\ee}{\end{equation}}
\newcommand{\bq}{\begin{equation*}}
\newcommand{\eq}{\end{equation*}}
\newcommand{\ba}{\begin{align*}}
\newcommand{\ea}{\end{align*}}
\newcommand{\bp}{\begin{proof}}
\newcommand{\ep}{\end{proof}}
\newcommand{\bL}{\begin{lemma}\label}
\newcommand{\eL}{\end{lemma}}
\newcommand{\bP}{\begin{proposition}\label}
\newcommand{\eP}{\end{proposition}}
\newcommand{\bC}{\begin{corollary}\label}
\newcommand{\eC}{\end{corollary}}
\newcommand{\bT}{\begin{theorem}\label}
\newcommand{\eT}{\end{theorem}}
\newcommand{\bR}{\begin{remark}\label}
\newcommand{\eR}{\end{remark}}
\newcommand{\bD}{\begin{definition}\label}
\newcommand{\eD}{\end{definition}}
\DeclareMathOperator{\tr}{Tr}
\DeclareMathOperator{\ra}{rank}

\newcommand{\A}{\mathscr{A}}
\newcommand{\Cu}{\mathscr{O}}
\newcommand{\Mul}{\mathscr M}
\title[ ]{Finite sums of projections in purely infinite\\ simple
C*- algebras  with torsion $K_0$}

\author{Victor Kaftal}

\address{Department of Mathematics\\
University of Cincinnati\\
P. O. Box 210025\\
Cincinnati, OH\\
45221-0025\\
USA}

\email{kaftalv@ucmail.uc.edu}

\author{P. W. Ng}

\address{Department of Mathematics\\
University of Louisiana\\
217 Maxim D. Doucet Hall\\
P.O. Box 41010\\
Lafayette, Louisiana\\
70504-1010\\
USA}

\email{png@louisiana.edu}

\author{Shuang Zhang}

\address{Department of Mathematics\\
University of Cincinnati\\
P.O. Box 210025\\
Cincinnati, OH\\
45221-0025\\
USA}

\email{zhangs@math.uc.edu}

\begin{document}

\begin{abstract}
Assume that  $\A$ is a purely infinite simple C*-algebra whose $K_0$
is a torsion group, namely, contains no free element. Then a
positive element $a\in \A$ can be written as a finite sum of
projections in $\A$ if and only if either $a$ is a projection or
$\|a\|>1$.
\end{abstract}

\maketitle

\section{\bf Introduction}
In \cite[Theorem 1.1]  {Fp69} Fillmore proved that a positive
finite-rank operator   $a\in B(H)$ can be expressed as finite sums of
projections  if and only if $\tr(a)\in \mathbb N$ and $\tr(a)\ge
\ra(a)$, where $\tr $ is the standard integer-valued trace.

An
immediate natural question is: which infinite-rank operators can
also be expressed as finite sums or infinite sums of
 projections,
 where the infinite sums of projections are assumed to converge in the strong topology.
  Motivated by needs in Frame Theory, this question was posed by
 Dykema, Freeman, Kornelson,  Larson, Ordower,  and Weber in   (\cite {DFKLOW}),
 where they also proved a sufficient condition.

Recently,   we obtained in \cite[Theorem 1.1] {KNZW*09} the  complete characterization: a positive operator $a\in B(H)$ is a
(possibly infinite) sum of projections if and only if $$\tr(a_+)\ge
\tr(a_-)\quad\text{and if}\quad   \tr(a_+)< \infty\quad\text{then}\quad   \tr(a_+)- \tr(a_-)\in \mathbb
N\cup\{0\}.$$ Here  $a_+:=(a-I)\chi_{(1, \infty)}(a)$,
and $a_-:= (I-a) \chi_{[0,1)}(a)$ where $\chi(a)$ is the spectral measure of $a$.

A similar characterization, but
without the integrality condition, can be given for all positive
diagonalizable operators in $\sigma$-finite von Neumann factors. In
particular, in type III factors  a positive operator $a$ can be
expressed as a sum of projections if and only if
  either $a$ is a projection or $\|a\|>1$.

Turning into another direction,  we generalized the  $B(H)$ results to  bounded module maps on Hilbert C*-modules. In  \cite{KNZC*09} we considered  the multiplier algebras of purely infinite and simple C*-algebras, equipped with the strict topology (in analogy of the strong operator topology on
$B(H)$)  and obtained the following complete characterization (\cite[Theorem 1.1 ]{KNZC*09}):  if $\A $ is a
$\sigma$-unital, nonunital, purely infinite simple  $C^*$-algebra
and $a$ be a positive element of $\Mul(\A )$, then $a$ is a  sum of
projections belonging to $\A$, with convergence in the strict
topology, if and only if one of the following mutually exclusive
conditions hold:
\begin{enumerate}
\item[(i)] $\| a\|_{e} > 1$
\item[(ii)]) $\| a \|_{e} = 1$ and $\|a\| > 1$.
\item[(iii)] $a\in \Mul (\A )\setminus \A $ is a projection.
\item[(iv)]  $a$ is the sum of finitely many projections belonging to $\A$.
\end{enumerate}
Here  $\| \cdot \|_{e}$ denotes the essential norm on $\Mul(\A )$, namely  $\| a \|_{e}=\|\pi(a)\|$, where
$\pi:\Mul(\A)\to \Mul(\A)/ \A$ is the canonical quotient map.

A harder problem  is to characterize  the case (iv) above, that is,
which positive elements in a C*-algebra $ \A$ are finite sums of
projections in $\A$. This question is still open even in $B(H)$.

A sufficient condition for $a \in B(H)^+$  to be a finite sum of
projections is that $\|a\|_e>1$. This result was reported in a survey article \cite
[Theorem 4.12] {Wpy94} by Wu quoting unpublished joint work of Choi
and Wu in 1988. The proof of this sufficient condition was presented  in their recent paper
\cite[Theorem 2.2] {WpCm09}.

The special but important case, $\alpha I$   with $\alpha
>1$,  follows also from a delicate analysis  in \cite {KRS02}, \cite {KRS03}, and
several other papers which characterized the set
$$ \Sigma_n:= \{\alpha  \in \mathbb R \mid   \alpha I \text{ is the sum of $n$ projections} \}$$
A  consequence of their analysis yields that $\bigcup_1^\infty
\Sigma_n= [1, \infty)$, i.e., $\alpha I$ is a finite sum of
projections if and only if $\alpha \ge 1$.

An independent and different proof of the sufficient condition  by
Choi and Wu was given by the authors of the present paper in
  \cite {KNSfinsumW*}: If $M$ is an infinite $\sigma-$finite von Neumann factor and $a\in M$ is positive, then
  a sufficient condition for $a$ to be a finite sum of projections in $M$ is, respectively,
\begin{enumerate}
\item [(i)] $\|a\|_e>1$ when $M$ is of type I$_\infty$ (where $\|\cdot \|_e$ denotes the usual essential norm of $B(H)$);
\item [(ii)] $\|a\|_e>1$ when $M$ is of type II$_\infty$  (where $\|\cdot
\|_e $ denotes the essential norm relative to the closed ideal
  $ J(M)$ generated by all finite projections of $M$);
\item [(iii)] $\|a\|>1$ when $M$ is of type III. Moreover,  $a\in M^+$ is a
finite sum of projections if and only if either $a$ is a projection
or $\|a\|>1$.
\end{enumerate}
By  employing  a suitable notion of essential central spectrum, we proved also
that a  similar sufficient condition holds  also for
global properly infinite von Neumann algebras  \cite[Theorem 4.3]
{KNSfinsumW*}.

For multiplier algebras we proved in \cite[Theorem 3.3] {KNZC*10}  a
similar result: if $\A $ is a $\sigma$-unital,
  nonunital, purely infinite simple C*-algebra, then every  positive element  $a\in \Mul(\A )\setminus \A $ with
$\|a\|_{e}>1$ can be expressed as  a finite sum of projections.

Since   $\sigma$-unital purely infinite simple C*-algebras   share
some properties of $\sigma$-finite type III factors, one would hope
that every positive element $a$   in a purely infinite simple
C*-algebra with $\|a\|>1$ could also be written as a finite sum of
projections. And if not, one would like to see which such algebras
do not have this property.

The objective of this paper is to prove  the following result.

\bT{T:main} Let $\A$ be a purely infinite simple C*-algebra such
that $K_0(\A )$ is a torsion group, namely, contains no free
elements, then  a positive element $a\in \A^+$ is a finite sum of
projections if and only if either $a$ is a projection or $\|a\|>1$.
\eT

\medskip
Both in the case of von Neumann algebras  (\cite{KNSfinsumW*}) and
in the case of multiplier algebras (\cite{KNZC*09} and
\cite{KNZC*10}), the results were obtained by first
decomposing $a$ into a special infinite sum of projections
(converging in the strong topology in the von Neumann algebra case,
and in the strict topology in the multiplier algebra case), and
then, reassembling these sums into finite sums of projections.
However, for  purely infinite simple C*-algebras  the above
techniques can no longer be applied, and we will have to use a different approach.

A key reduction is that if a purely infinite simple C*-algebra $\A$ has
the property that $\alpha p$ is a finite sum of projections  for all $\alpha >1$ and all
projections $p\in \A$, then  each   $a\in \A^+$ with $\|a\|>1$ is
also of a finite sum of projections. The key  in this  reduction is  the  fact derived in (\cite{KNZC*10}) that each positive element $a\in \A$ can be written as positive linear combination of projections.

Then we prove that in any purely infinite simple C*-algebra,  $\alpha p$ is a finite sum of projections for any $\alpha
>1$ and any projection $p$ that represents the zero element of $K_0(\A )$ (Proposition \ref
{P:real}).

This will be done by showing that $\alpha p$ is of a finite sum of projections
for any $\alpha \in [\frac{3}{2}, 3]$ (Lemma \ref {L:big}), and
then, prove  that $\alpha p\oplus \beta q$ is a finite sum of projections for any pair of
rational numbers $\alpha >1$ and  $\beta\ge 0$ and every projection $p\sim q$, which represents the zero element of  $K_0(\A )$ (Lemma \ref {L:rational}).

Combining these results, we then obtain that  $\alpha p\oplus \beta q$ is a finite sum of projections for any $\alpha >1$, $\beta\ge 0$,  and any  projection  $p, q\in \A$ which represent  elements of $K_0(\A )$ with finite orders (Proposition \ref {P:O_n}).

As a consequence, if  all elements of $K_0(\A )$ are of finite order, i.e.,$K_0(\A )$ is a torsion group, then
$a\in \A^+$ is a finite sum of projections  if and only if either $a$ is a projection
  or $\|a\|>1.$ (Theorem  \ref {T:main}).

\medskip

For the readers'  convenience we summarize some  basic facts on K-theory of
infinite simple C*-algebras that we will use in this paper.
 A handy reference  are the original  papers by  Cuntz \cite {Cj77} and \cite{Cj81}.\\
Recall  that   $[p]$ denotes the equivalence class of a projection
$p\in \A$ for any C*-algebras $\A $. If  $\A $ is infinite and
simple, then the $K_0$ group of $\A$ can be described as follows:
\begin{enumerate}
\item
$K_0(\A) = \{[p]\mid p \text{ is an infinite projection in }\A\}.$
\item
The addition `+' on $K_0(\A)$ is defined as follows: for any two infinite projections $p,q\in \A$, choose $p'\sim p$, $q'\sim q$ so
that $p'q'=0$. then $[p]+[q]=[p'+q']$. In particular, we denote  $$n[p]=\overset{ n \text{ times}} {\overbrace{[p]+[p]+\cdots [p]}}.$$
\item
The zero (or unit) element $\theta\in K_0(\A)$ is defined as
$[p-p']$ for any pair of projections $p'< p$ with $p'\sim p$ and
$p-p'$
 infinite.

\item
For any infinite projection $q\in \A$ and any positive integer $m$,
$q$ can be decomposed into the sum $q=\sum_{j=0}^mq_j$ of
 mutually orthogonal projections $q_j$, with $q_0\sim q$ and $[q_j]=\theta$ for all $j>0$.
\item An infinite C*-algebra $\A$ has torsion $K_0$ if for every projections $p\ne 0$ the
element $[p]$ is finitely cyclic in $K_0(\A)$,
 namely, there is a positive integer $n$ for which $n[p]=\theta$,  the zero element of
 $K_0(\A)$. Among algebras with torsion $K_o$ are the Cuntz
 algebras $\Cu_n$ for $2\le n< \infty$, for which  $K_0(\Cu_n)= \mathbb Z/(n-1)\mathbb Z$;
 in particular, $K_0(\Cu_2) =\{\theta\}$, the trivial group.
\end{enumerate}

\section{\bf Finite sums of projections and K-theory}

The whole section will be devoted to establish the proof of our main
result by a sequence of technical lemmas.

 \bL{L:big}
 Let $\A $ be an  infinite simple  C*-algebra,   $p\in \A$ be a projection that represents the zero
element of $K_0(\A )$  (i.e.,  $[p]=\theta$) and $\gamma\in [\frac{3}{2}, 3]$. Then  there are six projections $q_j\in \A$,
 with $q_j\sim p$ and such that
 $ \gamma p =\sum_1^{6}q_j.$
 \eL

 \bp
Decompose $p=p_1+p_2$ into the sum of two mutually orthogonal projections $p_1\sim p_2\sim p$. Let
\ba
a_1&:= (2\gamma- 3)p_1+(3-\gamma)p_2\\
a_2&:= (3-\gamma)p_1+(2\gamma- 3)p_2
\end{align*}
 Then
 $$\gamma p=a_1+a_2\qquad\text{and}\qquad a_1\ge 0, ~a_2\ge 0.$$
Further decompose $p_2=p_2' +p_2''$ into the sum  of two mutually orthogonal projections $p_2' \sim p_2''\sim p_2$. Then
$$a_1= (2\gamma- 3)p_1+ (3-\gamma)p_2'+(3-\gamma)p_2''.$$
Since the three projections $p_1$, $p_2'$, and $p_2''$ are mutually
orthogonal and equivalent, $a_1$ belongs to an embedding of $\mathbb
M_3(\mathbb C)$ in $p\A p$ where $p_1$, $p_2'$, and $p_2''$ are
identified  with rank-one projections. Thus $a_1$ can be identified
with a positive matrix in $\mathbb M_3(\mathbb C)$ with $\tr(a_1)=
3$ and $\ra(a_1) \le 3$. Then it follows from  \cite [Theorem 1.1] {Fp69} that $a_1$ is the sum of three
projections (all equivalent to $p$.) The same argument applies to
$a_2$, and hence,  the proof is completed.

\ep

\bL{L:rational} Let $\A $ be an infinite simple C*-algebra, $p,q\in
\A$ be two orthogonal projections that represent the zero
element of $K_0(\A )$,  and let $\alpha, \beta$ be rational numbers with $\alpha >1$  and  $\beta\ge 0$. Then $a:= \alpha p+\beta q$ is a finite sum of projections.
 \eL
 \bp
 Set $\alpha= \frac{k}{m}$ and $\beta= \frac{h}{m}$ with $k,h,m$ integers and choose an integer $r \ge \frac{m}{k-m}$,
  equivalently, $rk\ge rm+m$.
 Then decompose $p= \sum_1^{rm} p_j$ (resp., $q = \sum_1^{m} q_j$) into the sum of $rm$ (resp. $m$) mutually
 orthogonal projections
 $p_j\sim q_j\sim p$. Then
 $$a= \frac{k}{m}\sum_1^{rm} p_j+ \frac{h}{m}\sum_1^{m} q_j$$ can be identified with a positive matrix in
 $\mathbb M_{rm+m}(\mathbb C)_+$ with $$\tr(a)=rk+h\in \mathbb N\quad \text{and}\quad   \ra(a)\le rm+m.$$ Since $rk+h\ge rm+m$, it follows once again from \cite [Theorem 1.1] {Fp69}  that $a$ is a finite sum of
projections.

\ep

\bP{P:real} Let $\A $ be an infinite simple C*-algebra,  $p,q\in
\A$ be two orthogonal equivalent projections that represent the zero
element of $K_0(\A )$, and  let $\alpha, \beta$ be real numbers with $\alpha >1$  and  $\beta\ge 0$. Then $a:= \alpha p+\beta q$ is a finite sum of projections.
 \eP
 \bp
If $\beta=0$, write $\alpha p = \alpha p_1+\alpha p_2$,
where $p_1\sim p_2\sim p$ and $p_1p_2=0$. Thus assume without loss of generality that $\beta >0$.

 Let $\rho_1, \rho_2, \rho_3$ be rational numbers with
\ba 1&< \rho_1 < \alpha\\
\frac{1}{3}(\alpha-\rho_1)&< \rho _2 < \frac{2}{3}(\alpha-\rho_1)\\
\frac{1}{3}\beta&< \rho _3 < \frac{2}{3}\beta.
\end{align*}
Decompose $p= p'+p''$ into the sum of two mutually orthogonal projections $p'\sim p''\sim p.$ Then
$$a= \rho_1p +(\alpha-\rho_1)p'' +(\alpha-\rho_1)p' + \beta q.$$
Since $\frac{\alpha-\rho_1}{\rho_2}\in [\frac{3}{2}, 3]$  and
$\frac{\beta}{\rho_3}\in [\frac{3}{2}, 3]$, by Lemma \ref {L:big} we
can find projections $\{q_j\}_1^{18}$  with $q_j\sim p$ for all  $j$
and such that \ba
\frac{\alpha-\rho_1}{\rho_2}p''&= \sum_1^6q_j,\\
\frac{\alpha-\rho_1}{\rho_2}p'&= \sum_7^{12}q_j,\\
\frac{\beta}{\rho_3}q&= \sum_{13}^{18}q_j.
\end{align*}
It follows that
$$a = \rho_1p+ \rho_2 \sum_1^{6}q_j+  \rho_2\sum_7^{12}q_j+\rho_3 \sum_{13}^{18}q_j
$$
Now decompose $p'=\sum_1^{6}p_j$ and $p''=\sum_7^{18}p_j$  into the sum of mutually orthogonal projections $p_j\sim p$. Thus
$$
a=
 \sum_1^{6}\big( \rho_1p_j+ \rho_2q_j\big)+ \sum_7^{12}\big( \rho_1p_j+ \rho_2q_j\big) +\sum_{13}^{18}\big( \rho_1p_j+ \rho_3q_j\big).
$$
Since
all the coefficients $\rho_i$ are rational, $\rho_1>1$, and  $p_jq_j=0$ for $1\le j\le 18$ because $p'  p''= p q = 0$, it follows  by Lemma
\ref {L:rational} that each of the above $18$ summands is a finite sum of projections, and hence, so is $a$.

 \ep

 Taking $\beta = 0$ in the above proposition we thus have:
 \bC{C:K_0}
Let $\A $ be an infinite simple C*-algebra and $p$ be a projection
in $ \A$ that represents the zero element of $K_0(\A)$. Then   for
all $\alpha>1$ the positive element $\alpha p$ is a finite sum of
projections.
 \eC

This result extends also to projections that represent finitely cyclic elements of $K_0(\A).$

 \bP{P:O_n} Let $\A$ be a purely infinite simple C*-algebra
 and  $p $ be a projection of $\A$ that represents a finite order element  of $K_0(\A)$ (i.e.,  such that $[p]$ is finitely cyclic). Then $\alpha p$ is a finite sum of  projections.
 \eP
\bp Assume also without loss of generality that $\alpha < 2$.
Let $\theta$ be the zero element of $K_0(\A)$.
 By hypothesis  $n[p]=\theta$ for some integer $n\in \mathbb N$.
 Then $\frac{1}{\alpha -1}- \frac{1}{n}>0$ and let $m\ge 0$ be the integer part of $\frac{1}{\alpha -1}- \frac{1}{n}$  and
 $0\le \delta < 1$ be its fractional part, i.e.,
$$
m = \frac{1}{\alpha -1}- \frac{1}{n} -\delta.
$$
Since $[p]= (m+1)\theta +[p] = ((m+1)n+1)[p]$, we can decompose $p$ into
the sum $ p= \sum_{j=0}^{(m+1)(n)}p_j$ of mutually orthogonal
projection $p_j\sim p$. Then
$$\alpha p = \alpha\sum_{j=0}^  {mn} p_j + (\alpha-1)\delta \sum_{j=
 mn+1}^{(m+1)n} p_j + \big(\alpha - (\alpha-1)\delta\big)\sum_{j=  mn+1}^{(m+1)n} p_j .$$
Let
$$b:=  \alpha\sum_{j=0}^  {mn} p_j + (\alpha-1)\delta \sum_{j=  mn+1}^{(m+1)n} p_j $$
and
$$
c:= \big(\alpha - (\alpha-1)\delta\big)\sum_{j=  mn+1}^{(m+1)n} p_j .
$$
Then $a=b+c$. Since all the projections $p_j$ are mutually
orthogonal and equivalent (to $p$), the positive operator  $b$ can
be identified with a finite matrix with $\ra(b) \le (m+1)n+1$ and
with
$$\tr(b) =  (mn+1)\alpha + n(\alpha-1)\delta= (m+1)n+1.$$ Thus, by
\cite [Theorem 1.1] {Fp69}, $b$ is a finite sum of
projections. Since $\alpha - (\alpha-1)\delta>1$ and since
$$[\sum_{j= mn+1}^{(m+1)n} p_j]=n[p]=\theta,$$ by Corollary \ref {C:K_0} $c$ too is a finite sum of projections,
and hence, so is $a$. \ep

In particular, if  $\A $ is a purely infinite simple
C*-algebra such that $K_0(\A )$ is a torsion group, i.e.,  every
element is of finite order, then $\A$ has the property that the  element $\alpha p$ is a  finite sum of projections for every $\alpha >1$ and every projection $p\in \A$. As we will see, this property extends to all positive operators with norm larger than 1.

\bL {L:extension} Let $\A$ be a purely infinite simple C*-algebra
with the property that the  element $\alpha p$ is a  finite sum of projections for every $\alpha >1$ and every projection
$p\in \A$.
Then for every orthogonal pair of projections $0\ne p, q\in \A$ and
every pair of real numbers $\alpha, \beta$ with $\alpha>1$ and $\beta \ge 0$, the operator $a= \alpha p +
\beta q$ is a  finite sum of projections. \eL

\bp Without loss of generality assume that $\beta <1$.
Let $n\in
\mathbb N$ be the integer part of $\frac{2-\beta}{\alpha -1}$ and $0\le \delta <1$ be its fractional part and
let $\epsilon: =\delta(\alpha-1)$. Then $0\le  \epsilon< \alpha -1$ and $$n\alpha
+\beta +\epsilon = n+2.$$ Since $[p]=(n+1)[q]+\big([p]-(n+1)[q]\big)$ in
$K_0(\A)$, we can decompose $p$ into the sum
 $p=p'+ \sum_{j=0}^nq_j$ of mutually orthogonal nonzero projections with $q_j\sim q$ for $0\le j\le n$ and a nonzero ``remainder" $p'$ also orthogonal to all the projections $q_j$. Then
$$
a=  \alpha p' + (\alpha-\epsilon)q_0+   \epsilon q_0 +
 \sum_1^n \alpha q_j+\beta q
$$
Let $b:= \epsilon q_0+  \sum_1^n \alpha q_j+\beta q$. Then $b\ge 0$
and since all the projections are mutually orthogonal and
equivalent, we identify it with a positive matrix with
$$\tr(b)=n\alpha +\beta +\epsilon = n+2\ge \ra(b).$$ Thus $b$ is a
 finite sum of projections. By hypothesis so are  $\alpha p'$ and $
(\alpha-\epsilon)q_0$ and hence so is $a$.

\ep

 \bT{T:general} Let $\A$ be a purely infinite simple C*-algebra with
 the property that the  element $\alpha p$ is a  finite sum of projections for every $\alpha >1$ and every projection
$p\in \A$.
 Then
$a\in \A^+$ is a  finite sum of projections if and only if either
$a$ is a projection or $\|a\|>1$. \eT
\bp If $a$ is a sum of
projections, then  $\|a\|\ge 1$. Then either the projections are
mutually orthogonal, in which case $a$ itself is a projection, or
$\|a\| > 1$.

For the converse implication, assume that $\|a\| > 1$. By passing if necessary to the hereditary
algebra $(a\A a)^-$ generated by $a$ which is $\sigma$-unital, we assume without loss of generality
 that $\A$ itself is $\sigma$-unital. By \cite [Remark 2.1 (II)]{KNZC*10}, $\A$ has the PLP property,
 namely the positive cone of each of its hereditary C*-subalgebras is the closure of the positive combinations
of its projections. But then, by \cite[Lemma 2.10]{KNZC*10}, either
$a= \alpha p\oplus  c $ for some $\alpha>1$ and nonzero projections
$p \in \A$ and a positive element $c\in\A$, or $a$ is the sum of two
such operators. Assume without loss of generality that  $a= \alpha
p\oplus  c $.
 By \cite [Theorem 2.11]{KNZC*10}, $c=\sum_1^n\lambda_jq_j$ for some scalars $\lambda_j>0$ and projections $q_j.$
  Decompose $p$ into the sum $p= \sum_1^np_j$ of mutually orthogonal nonzero projections $p_j\in \A$. Then
$$
a= \sum_1^n\big(\alpha p_j+ \lambda_jq_j)
$$
where $p_jq_j=0$ for all $j$. By Lemma \ref {L:extension}, each
summand $\alpha p_j+ \lambda_jq_j$ is a  finite sum of projections,
and hence, so is $a$. 

\ep

Combining Proposition  \ref{P:O_n} and Theorem \ref {T:general}  we thus obtain our main result.

\

\noindent \bf{Theorem \ref{T:main}} \it{Let $\A$ be a purely infinite simple C*-algebra such
that $K_0(\A )$ is a torsion group, namely, contains no free
elements, then  a positive element $a\in \A^+$ is a finite sum of
projections if and only if either $a$ is a projection or $\|a\|>1$.
}

\bR{R:W*} \item [(i)] Purely infinite
simple C*-algebras with $K_0(\A )$ a torsion group include of course all those trivial $K_0$ and  all Cuntz algebras $\Cu_n$
with $2\le n< \infty$.
\item [(ii)] The  proofs of Lemmas \ref {L:big}, \ref {L:rational}, and
Proposition \ref {P:real},  and hence, the conclusion in Corollary
\ref {C:K_0}, hold also in $\sigma$-finite von  Neumann factors,
 if restricted to infinite projections. Namely,   if $M$ is $\sigma$-finite von  Neumann factor,
then $\alpha p$ is a  finite sum of projections for all $\alpha>1$
and all infinite projections $p$. This observation provides a simple alternative proof for the fact that $\alpha 1$
is a finite sum of projections in $B(H)$ if $\alpha
>1$ (see \cite {KRS02}, \cite {KRS03}, \cite {Wpy94}, \cite
{WpCm09}, \cite {KNSfinsumW*}.)
\item [(iii)] Notice explicitly that the techniques of the present paper  do not provide an alternative proof for the sufficiency of the condition $\|a\|_e>1$ for $a\in B(H)^+$ to be a  finite sum of projections. Indeed, while by a  proof similar to the one for Theorem \ref {T:main} (see \cite [Theorem 4.3]{KNSfinsumW*})
it is enough to consider elements of the form  $a=\alpha p\oplus \beta q$,  where $\alpha>1$, $\beta \ge 0$, $p, q$ are projections,
and $p$ is infinite, it is not enough to consider only elements where also $q$ is infinite. But the proofs of Proposition \ref {P:real} or
of Lemma \ref {L:extension} do not apply unless
  $q$ too is infinite;  a different technique is required when $q$ is finite (see  \cite [Lemma 4.2]{KNSfinsumW*}).

Notice  also  that the condition $\|a\|_e>1$ implies that $\pi(a)$ is a  finite sum of projections in the Calkin
 algebra $B(H)/K(H)$, but this alone does not imply that the same holds for $a$ (e.g., there are cases of
 operators $1+k$ with $k\in K(H)^+$ that are not  finite sum of projections (see \cite [Example 5.9]{KNSfinsumW*}).
\item [(iv)] Technical difficulties occur at attempting the analysis of purely infinite simple C*-algebras whose $K_0$ group contains free elements;
in particular, even $\Cu _\infty $ resists all efforts so far.
 Unlike all type III von Neumann factors in which all hereditary W*-subalgebras
 are unital,  hereditary C*-subalgebras of a $\sigma$-unital purely infinite simple C*-algebra $\A $ are
  either unital or stable (namely, $\A \cong \A \otimes K(H)$), where  $\A $ can be chosen unital and
  $K(H)$  is the stable algebra
  of all compact operators on a separable Hilbert space  (see \cite{Zhang1} and \cite{Zhang2}).
  It is well known that infinite  rank positive operators  in  $K(H)$  are not finite sums of
projections, but the same positive elements in the  canonical embedding of $K(H)$ as a C*-subalgebra of $\A \otimes K(H)$, become finite sums
of projections as long as $K_0(\A )$ is a torsion group. Whether the
same conclusion can be
    expected or not, when $K_0(\A )$ is not a torsion
    group, remains mysterious. If not, these purely infinite simple C*-algebras would exhibit quite a different behavior from  type III factors concerning the decompositions of its positive elements.
\eR

\

 {\it Acknowledgement: S. Zhang was supported by a Taft Center Travel Grant when the article was presented in
 Beijing, China in the summer of 2010.}

\end{document}